\documentclass[reqno]{amsart}

\usepackage{amssymb}
\usepackage{fullpage}
\usepackage{amsthm}
\usepackage{bm}
\usepackage{latexsym}
\usepackage{float}
\restylefloat{table}
\usepackage{amsmath}
\usepackage{multirow}
\usepackage{eufrak}
\usepackage{mathrsfs}
\usepackage[font={small,it}]{caption}
\usepackage{amscd}
\usepackage[all,cmtip]{xy}
\usepackage[usenames]{color}
\usepackage{graphicx}
\usepackage{color}
\usepackage{amscd}
\usepackage{float}
\usepackage{graphics}
\usepackage{tikz}
\usepackage{tikz-cd}
\usepackage{comment}
\usepackage{xspace}
\usepackage{mathtools}
\usepackage{lipsum}
\usepackage{tabularx}
\usepackage{booktabs}
\usepackage{amssymb}
\usepackage{amsthm}
\usepackage{graphicx}
\usepackage{subfig}
\usepackage{enumerate}

\usepackage[style=alphabetic, backend=bibtex, backref=true,
 doi=false,isbn=false,url=true,
 eprint=true, sorting=nyt, maxnames=99]{biblatex}
\usepackage{hyperref}
\usepackage[capitalize]{cleveref}

\usetikzlibrary{patterns,decorations.pathreplacing}
\usepackage{marginnote}

\newtheorem{theorem}{Theorem}[section]
\newtheorem{lemma}[theorem]{Lemma}

\theoremstyle{definition}
\newtheorem{definition}[theorem]{Definition}
\newtheorem{example}[theorem]{Example}

\newtheorem{prop}[theorem]{Proposition}

\theoremstyle{remark}
\newtheorem{remark}[theorem]{Remark}

\numberwithin{equation}{section}

\def\Gal{\mathrm{Gal}}
\def\Jac{\mathrm{Jac}}

\DeclareMathOperator{\Hom}{Hom}

\newcommand{\MA}{\mathcal{A}}
\newcommand{\MM}{\mathcal{M}}

\newcommand{\F}{\mathbb{F}}
\newcommand{\G}{\mathbb{G}}

\newcommand{\Z}{\mathbb{Z}}

\newcommand{\calK}{\mathcal{K}}

\newcommand{\calI}{\mathcal{I}}

\newcommand{\calH}{\mathcal{H}}

\newcommand{\vv}{\mathbf{v}}
\newcommand{\bbP}{\mathbf{P}}
\newcommand{\im}{\mathrm{im}}

\DeclareFontFamily{U}{wncy}{}
\DeclareFontShape{U}{wncy}{m}{n}{<->wncyr10}{}
\DeclareSymbolFont{mcy}{U}{wncy}{m}{n}
\DeclareMathSymbol{\Sh}{\mathord}{mcy}{"58}

\DeclarePairedDelimiter\floor{\lfloor}{\rfloor}
\DeclarePairedDelimiter\ceil{\lceil}{\rceil}

\newcommand{\htC}{\operatorname{ht}_{\mathrm{Car}}}

\newcommand{\hC}{h_{\mathrm{Car}}}
\newcommand{\perpC}{\perp_{\mathrm{Car}}}
\newcommand{\Span}{\operatorname{Span}}
\newcommand{\rank}{\operatorname{rank}}
\newcommand{\codim}{\operatorname{codim}}
\theoremstyle{plain}
\newtheorem{corollary}[theorem]{Corollary}
\theoremstyle{remark}
\crefname{prop}{proposition}{propositions}
\Crefname{prop}{Proposition}{Propositions}
\hypersetup{
  colorlinks=false,
  pdfborder={0 0 1},
  linkbordercolor={1 0 0},
  citebordercolor={1 0 0},
  urlbordercolor={1 0 0},
  bookmarksnumbered=true,
  pdftitle={A sharp bound for the a-number of a cyclic cover of the projective line},
  pdfauthor={Shiva Chidambaram}
}

\newcommand{\calG}{\mathcal{G}}

\title[A sharp bound for the $a$-number of a cyclic cover of $\bbP^1$]
{A sharp bound for the $a$-number of a cyclic cover of $\bbP^1$}

\author{Shiva Chidambaram}

\address{University of Wisconsin-Madison,
480 Lincoln Dr, Madison, WI 53706 USA.}
\email{chidambaram3@wisc.edu}
\urladdr{\url{https://people.math.wisc.edu/~chidambaram3/}}


\date{}

\begin{document}

\begin{abstract}
Let $X$ be a smooth projective curve of genus $g$ over an algebraically
closed field $k$ of characteristic $p>2$. We give upper bounds for the
$a$-number $a(X)$ when $X$ is hyperelliptic or, more generally, a tame cyclic
cover of the projective line. For hyperelliptic curves, writing
$g=pm+r$ with $0\le r<p$ and $n=(p-1)/2$, we prove the sharp bound
$a(X)\le nm+\min\{n,r\}$. For tame cyclic covers, we refine previously
known upper bounds and give curves attaining the resulting bound
in infinitely many genera in characteristics $3$ and $5$.
The proofs combine an alternating pairing on the bundle of locally exact 
differentials, originally due to Raynaud, with an explicit polynomial 
description of the Cartier operator. 
We also produce an explicit locus of hyperelliptic curves $X$ with 
$a(X)\ge a_g$, where $\frac{a_g}{g}$ converges to $\frac{p-1}{2p}$
as $g$ goes to $\infty$, whose dimension grows linearly with the 
genus $g$. These families exhibit excess intersections of the Torelli 
locus with loci of large $a$-number in the moduli space of principally 
polarized abelian varieties.
\end{abstract}
\maketitle

\section{Introduction}\label{sec:intro}

Let $p>2$ be a prime, and let $k$ be an algebraically closed field of
characteristic $p$. Let $X$ be a smooth projective
curve of genus $g$ over $k$. We write $a(X)$ for the $a$-number of
its Jacobian. It is an invariant of the $p$-torsion group scheme $\Jac(X)[p]$,
defined precisely as $a(X)=\dim_k\Hom_{k\text{-grp}}(\alpha_p,\Jac(X)[p])$, 
where $\alpha_p$ is 
the kernel of the Frobenius endomorphism on the group scheme $\G_a$.
Equivalently, by the Dieudonn\'e--de Rham comparison
\cite[Corollary~5.11]{Oda1969}, $a(X)$ is the dimension of
the kernel of Cartier on the space $H^0(X,\Omega_X^1)$ of regular 
differentials on $X$. The latter description allows us to study $a(X)$ 
using an equation for the curve. 

For hyperelliptic curves, the Cartier operator is particularly
explicit: its matrix is read off from the coefficients of a power
of the defining polynomial. The classical references are Cartier
and Yui \cite{Cartier1957,Yui1978}; we follow the semilinear
conventions explained by Achter--Howe \cite{AchterHowe2019}.
Garton--Thunder--Weir use function-field heights to determine the
distribution of $a$-numbers of hyperelliptic curves in characteristic
$3$ and the dimensions of their $a$-number loci
\cite{GartonThunderWeir2026}. For tame cyclic covers of $\bbP^1$,
Elkin decomposes the regular differentials into eigenspaces for the
Galois action and bounds the rank of Cartier on each eigenspace in
terms of the ramification data \cite[Theorem~1.1 and
Corollary~4.2]{Elkin2011}.

We obtain additional restrictions on the regular differentials killed
by Cartier. The geometric input is Raynaud's pairing on the bundle
of locally exact differentials \cite[\S4.1]{Raynaud1982}; see also
Tong \cite[\S\S1.2.1,2.2]{Tong2010}. On a hyperelliptic curve, the
involution acts by $-1$ on every regular differential. The pairing
of two regular Cartier-kernel differentials is therefore invariant, and hence vanishes
since it is again a regular differential. This observation bounds the 
number of independent rational directions represented by regular 
differentials in the Cartier kernel.

Let $\mathscr B_X$ denote Raynaud's bundle of locally exact
differentials on the Frobenius twist of $X$, and let $\rho(X)$ be
the rank of the saturation of its evaluation image. We recall these
constructions precisely in \cref{sec:orthogonality}. Put $n=(p-1)/2$.

\begin{theorem}\label{thm:main}
Let $X$ be a smooth hyperelliptic curve of genus $g\ge2$ over $k$.
Write $g=pm+r$, where $0\le r<p$. Then
\begin{equation}\label{eq:main-bounds}
 \rho(X)\le n,\qquad a(X)\le nm+\min\{n,r\}.
\end{equation}
For each fixed odd prime $p$, the rank bound is attained by curves
defined over $\F_p$ in every genus $g\ge p$. Both bounds are attained
simultaneously for infinitely many genera.
\end{theorem}

\begin{remark}
When $p=3$, \cref{thm:main} recovers the result $a(X)\le \ceil{g/3}$ of 
Garton--Thunder--Weir \cite[Theorem~1.1]{GartonThunderWeir2026}.
\end{remark}

\cref{thm:main} also recovers the following result of Ekedahl 
about superspecial hyperelliptic Jacobians.

\begin{samepage}
\begin{corollary}[{\cite[Theorem~1.1(ii)]{Ekedahl1987}}]
\label[corollary]{cor:ekedahl}
Let $X$ be a superspecial hyperelliptic curve of genus $g$ in characteristic $p>2$,
i.e., $a(X)=g$.
Then $g \le (p-1)/2$.
\end{corollary}

\begin{proof}
If $g>n$, the upper bound for $a(X)$ in \cref{thm:main} is strictly
less than $g$.
\end{proof}
\end{samepage}

Sharpness in \cref{thm:main} is supplied by the \emph{binomial family}
\begin{equation}\label{eq:intro-binomial}
 X\colon y^2=f(x)=A(x^p)+x^sB(x^p),\qquad 1\le s<p,
\end{equation}
where $A,B\in k[T]$ are nonzero and the right-hand side is squarefree
of degree $2g+1$ or $2g+2$. 
The name reflects the role of the binomial theorem in proving
that $\rho(X)=n$ for these curves in sufficiently large genus
(see \cref{subsec:binomial});
$A$ and $B$ need not be monomials. 
For $g\ge2$, let $\MM_g$ denote the moduli space of smooth
projective curves of genus $g$, and let
$\calH_g\subseteq\MM_g$ be the hyperelliptic locus.
For a non-negative integer $a$, the superscript $\ge a$
denotes the locus where the Jacobian has $a$-number at least $a$. 
For the statement below, let $\Lambda_s\subseteq\{0,\ldots,p-1\}$
be the set of least nonnegative residues of $0,s,\ldots,ns$
modulo $p$.
Then:

\begin{theorem}\label{thm:binomial-intro}
For a fixed odd prime $p$ and sufficiently large $g$, every smooth member
of \eqref{eq:intro-binomial} satisfies
\begin{equation}\label{eq:intro-binomial-a}
 \rho(X)=n,\qquad
 a(X)=a_{p,g,s}:=\sum_{\substack{0\le t<p\\t\notin\Lambda_s}}\floor*{\frac{g+t}{p}}.
\end{equation}
For every sufficiently large $g$, one can choose $s$ so that
$\calH_g^{\ge a_{p,g,s}}$ is nonempty and
\begin{equation}\label{eq:intro-codim}
 \codim_{\calH_g}\calH_g^{\ge a_{p,g,s}}
 \le\left(2-\frac4p\right)g+O_p(1),
 \qquad a_{p,g,s}=\frac{p-1}{2p}g+O_p(1).
\end{equation}
The corresponding locus in $\MM_g$ satisfies
\begin{equation}\label{eq:intro-Mg}
 \codim_{\MM_g}\MM_g^{\ge a_{p,g,s}}
 \le\left(3-\frac4p\right)g+O_p(1).
\end{equation}
\end{theorem}

For these thresholds $a_g=a_{p,g,s}$, the ratio $a_g/g$ tends to the
positive constant $(p-1)/2p$ as $g$ goes to infinity. Nevertheless, 
the locus $\MM_g^{\ge a_g}$ has codimension at most $(3-4/p)g+O_p(1)$ 
in $\MM_g$, or equivalently, dimension at least $4g/p-O_p(1)$.
This is strikingly in contrast with the locus
$\MA_g^{\ge a}$ in the moduli space $\MA_g$ of principally polarized
abelian varieties, which has pure codimension $a(a+1)/2$
\cite[Theorem~1.2 and \S11]{Oort2001} (see also \cite[\S12.3]{EkedahlVanderGeer2009}). 
The expected dimension of the intersection of the Torelli locus
with $\MA_g^{\ge a_g}$ is $3g-3-\frac{a_g(a_g+1)}{2}$,
which is negative for sufficiently large $g$. The binomial family 
nevertheless produces excess intersections, whose 
dimensions in fact grow at least linearly with the genus.
In characteristic $3$, Garton--Thunder--Weir
\cite[Theorem~1.1]{GartonThunderWeir2026} give the more precise
result that $\calH_g^{\ge a}$ has codimension $2a-1$ in
$\calH_g$ for $1\le a \le \ceil{g/3}$.

\begin{remark}
\label{rem:specialfamilies_P1xP1}
For fixed $\gamma\ge3$, similar sparse equations on
$\bbP^1\times\bbP^1$ give families of curves of gonality $\gamma$
in infinitely many genera, with $a(X)\ge a_g$ and
$a_g/g$ converging to a positive constant $\lambda_{p,\gamma}$ 
along this sequence of genera.
For sufficiently large genera in the sequence, the dimension of 
the corresponding locus $\MM_g^{\ge a_g}$ is at least a positive 
constant times $g$, where the constant depends only on $p$ and $\gamma$.
\end{remark}

We next consider a more general setting. Let $X$ be a connected tame 
cyclic cover of $\bbP^1$, with $\pi:X\to\bbP^1$ of degree $N\ge2$, and
$p\nmid N$.
Let $G=\Gal(X/\bbP^1) \cong \Z/N\Z$ and fix a faithful character
$\chi:G\to k^\times$. For $1\le i<N$, let $D_i$ be the
$\chi^i$-eigenspace in $H^0(X,\Omega_X^1)$, and put
$d_i=\dim_kD_i$ and $\mathbf d=(d_1,\ldots,d_{N-1})$.
Let $\sigma(i)$ be the integer in $\{1,\ldots,N-1\}$ congruent to
$p^{-1}i$ modulo $N$. The Cartier action $\mathcal C_X$ maps $D_i$ to $D_{\sigma(i)}$.
The pairing of regular exact differentials in the eigenspaces for
$\chi^i$ and its inverse $\chi^{N-i}$ is invariant, and hence zero.
This gives rank inequalities that can be used together with
Elkin's bounds for the Cartier rank.

Elkin's lower bound for the rank of Cartier on $D_i$
\cite[Corollary~4.2]{Elkin2011} implies that
$\dim_k\ker(\mathcal C_X|_{D_i})\le e_i$, where
\begin{equation}\label{eq:elkin-caps}
 e_i=d_i-\min\{2\lfloor d_i/p\rfloor,d_{\sigma(i)}\}.
\end{equation}
Summing over the eigenspaces gives \emph{Elkin's bound}
$a(X) \leq \mathsf E_p(\mathbf d):=\sum_{i=1}^{N-1}e_i$ 
\cite[Theorem~1.1]{Elkin2011}.

To describe the \emph{pairing bound} 
we first introduce the polynomial counting 
function $\beta_p(d,s)=sq+\min\{s,r\}$, where $d=pq+r$,
$0\le r<p$, and $s$ is an integer with $0\le s\le p$.
Writing a polynomial of degree less than $d$ as
$\sum_{j=0}^{p-1}v_j(x)^p x^j$ gives $r$ coordinate spaces of
dimension $q+1$ and $p-r$ of dimension $q$. The integer $s$
specifies how many coordinates are retained, and $\beta_p(d,s)$
is the largest possible sum of their dimensions. In particular,
$\beta_p(d,0)=0$ and $\beta_p(d,p)=d$.

If $s_i$ is the rank of the subbundle generated by the regular
Cartier-kernel differentials in $D_i$, the polynomial counting
argument gives $\dim_k\ker(\mathcal C_X|_{D_i})\le\beta_p(d_i,s_i)$.
Since regular Cartier-kernel differentials in the two eigenspaces $D_i, D_{N-i}$ corresponding to inverse characters 
pair trivially under Raynaud's pairing, we also have $s_i+s_{N-i}\le p-1$. We therefore consider 
the finite set
\[
 \mathscr R_{p,N}
 =\left\{\mathbf s\in\{0,\ldots,p-1\}^{N-1}:
          s_i+s_{N-i}\le p-1\text{ for every }1\le i<N\right\}.
\]
Maximizing the polynomial counting estimate over these integer
ranks gives the \emph{pairing bound}
\begin{equation}\label{eq:pairing-bound}
 \mathsf P_p(\mathbf d)
 =\max_{\mathbf s\in\mathscr R_{p,N}}
       \sum_{i=1}^{N-1}\beta_p(d_i,s_i).
\end{equation}

The two bounds can be imposed simultaneously: for each
eigenspace, we take the smaller of the pairing-based polynomial counting
bound and Elkin's kernel bound, before maximizing over the
ranks. This defines the \emph{combined bound}
\begin{equation}\label{eq:combined-bound}
 \mathsf B_p(\mathbf d)
 =\max_{\mathbf s\in\mathscr R_{p,N}}
       \sum_{i=1}^{N-1}\min\{\beta_p(d_i,s_i),e_i\}.
\end{equation}

\begin{theorem}\label{thm:cyclic-intro}
Let $X\to\bbP^1$ be a connected tame cyclic cover of degree $N$
over $k$, and let $\mathbf d=(\dim_kD_i)_{i=1}^{N-1}$ be the tuple
of dimensions of the $\chi^i$-eigenspaces of $H^0(X,\Omega^1_X)$. 
Then $a(X)\le\mathsf B_p(\mathbf d)$.
\end{theorem}

By construction, the combined bound $\mathsf B_p(\mathbf d)$
is no larger than either the pairing bound or Elkin's bound:
$\mathsf B_p(\mathbf d)
 \le\min\{\mathsf P_p(\mathbf d),\mathsf E_p(\mathbf d)\}$.
The examples below show that it can give a strict
improvement over both.

\begin{example}
In characteristic $3$, the smooth projective models
$X_t\colon y^5=x^{15t+5}-1$ with $t\ge1$
have genus $30t+6$ and $a$-number
$a(X_t)=\mathsf B_3(3t,6t+1,9t+2,12t+3)=11t+3$.
The pairing bound $\mathsf P_3$ and Elkin's bound $\mathsf E_3$
are $14t+4$ and $13t+4$, respectively. We compute these values in
\cref{subsec:cyclic-sharpness}.
\end{example}
\noindent This example proves sharpness of the combined bound for the ramification 
types it realizes. We do not assert optimality for every ramification 
type. \Cref{subsec:cyclic-sharpness} also gives curves attaining the
combined bound in infinitely many genera in characteristic $5$.

We also show that the pairing bound $\mathsf P_p$ itself improves 
on Elkin's bound $\mathsf E_p$ in the following ranges.

\begin{theorem}\label{thm:comparison-intro}
\begin{enumerate}[(i)]
\item If $2\le N<p$, every connected tame cyclic cover
$X\to\bbP^1$ of degree $N$ over $k$ satisfies
$\mathsf P_p(\mathbf d)\le\mathsf E_p(\mathbf d)$.
\item Fix $p\ge5$ and an integer $N\ge2$ with $p\nmid N$.
There is an integer $g_0(p,N)$ such that every connected tame
cyclic cover $X\to\bbP^1$ of degree $N$ and genus $g\ge g_0(p,N)$
satisfies $\mathsf P_p(\mathbf d)<\mathsf E_p(\mathbf d)$.
\end{enumerate}
\end{theorem}
\noindent Part (i) restricts the covering degree but allows every genus.
Part (ii) allows any fixed covering degree prime to $p$, but assumes
$p\ge5$ and sufficiently large genus. 

Here is a brief outline of the article. \Cref{sec:orthogonality} introduces the 
alternating pairing on the bundle of
locally exact differentials and proves the fundamental geometric rank inequalities.
\Cref{sec:polynomial} relates these ranks to polynomial syzygies and
gives an exact formula for $a(X)$ in terms of the successive minima of an explicit 
polynomial lattice, for a hyperelliptic curve $X$. 
\cref{subsec:coordinate-pairing} also gives an explicit polynomial coordinate 
description of Raynaud's alternating pairing in the hyperelliptic case.
\Cref{sec:cyclic} proves the upper bounds for tame cyclic covers
and compares them with Elkin's bound. The hyperelliptic binomial
family and the sharpness examples are treated in \cref{sec:families}.

\subsection*{Notation and conventions.}
All curves are smooth, projective, and geometrically connected.
Smooth projective models are described by affine equations. We use
$X$ for a curve, $x,y$ for coordinates, and $T$ as a formal
polynomial variable. The notation $a(X)$ always denotes the
$a$-number of $\Jac(X)$; a letter $a$ in a locus $\MM_g^{\ge a}$
is a numerical threshold. All moduli statements concern genus $g\ge2$.
All asymptotic constants may depend on the fixed characteristic,
covering degree, or gonality. All maxima over rank parameters are
over integers. We put $\deg0=-\infty$ and $k[x]_{\le d}=0$ for $d<0$.

\subsection*{Acknowledgments}
I thank Jordan Ellenberg for our discussions and his questions 
which have been greatly useful for this work. 
I also thank Dima Arinkin and Joshua Mundinger for useful conversations.
I also gratefully acknowledge partial support by NSF grant DMS-2301386.

\subsection*{AI use acknowledgment}
Through the ChatGPT for Academic Researchers plan, GPT-6 Pro was used 
to assist with literature searches, and computational checks. 
It was also used to assist in drafting the manuscript, by building on 
a rough working draft written by the author. Extensive rewrites and revisions were done
afterwards, both directly by the author, and by specifically directing 
AI to do them, which were then verified by the author.
The author takes full responsibility for the content of the article.

\section{Orthogonality of regular exact differentials}
\label{sec:orthogonality}

\subsection{Locally exact differentials and global sections}

Let $X_1=X^{(p)}$ be the Frobenius twist of $X$, and let
$F_X:X\to X_1$ be relative Frobenius. We write $F_X^\#$ for the
induced map on functions. The bundle $\mathscr B_X$ of locally exact
differentials and the relative Cartier operator $c_X$ fit into
Raynaud's exact sequences \cite[\S4.1]{Raynaud1982}
\begin{align}\label{eq:cartier-sequences}
 \begin{split}
 0&\longrightarrow\mathcal O_{X_1}
  \xrightarrow{F_X^\#}F_{X*}\mathcal O_X
  \xrightarrow{d}\mathscr B_X\longrightarrow0,\\
 0&\longrightarrow\mathscr B_X
  \longrightarrow F_{X*}\Omega^1_{X/k}
  \xrightarrow{c_X}\Omega^1_{X_1/k}\longrightarrow0.
 \end{split}
\end{align}
Thus $\mathscr B_X$ has rank $p-1$ on $X_1$. Both the source and
target of $c_X$ are sheaves on $X_1$. 
We reserve $\mathcal C_X$ for the Cartier operator acting on
differentials on $X$, obtained from $c_X$ by identifying the
differentials on $X_1$ with those on $X$ through inverse Frobenius
on $k$. The resulting operator is $p^{-1}$-semilinear:
$\mathcal C_X(\lambda\omega)=\lambda^{1/p}\mathcal C_X(\omega)$
for $\lambda\in k$.

Let $\alpha_p=\ker(F:\mathbf G_a\to\mathbf G_a)$. As mentioned in 
the introduction, the $a$-number of $X$ is defined by the first 
expression below, while the second follows from the Dieudonn\'e--de Rham 
comparison \cite[Corollary~5.11]{Oda1969}:
\begin{align}\label{eq:a-definitions}
 a(X)=\dim_k\Hom_{k\text{-grp}}(\alpha_p,\Jac(X)[p])
     =\dim_k\ker(\mathcal C_X|_{H^0(X,\Omega_X^1)})
     =h^0(X_1,\mathscr B_X).
\end{align}
The last equality follows from the second exact sequence in
\eqref{eq:cartier-sequences}. A global section of $\mathscr B_X$
is a globally regular differential on $X$ with regular primitives
locally. It need not have a globally regular primitive.

Let $\calI_X\subseteq\mathscr B_X$ be the subsheaf generated by the
global sections, namely the image of the evaluation map:
\[
 \calI_X=\im\left[H^0(X_1,\mathscr B_X)\otimes_k\mathcal O_{X_1}
       \longrightarrow\mathscr B_X\right].
\]
Define $\calG_X$ to be the saturation of the evaluation image
$\calI_X$ in $\mathscr B_X$, and put $\rho(X)=\rank\calG_X$.
Explicitly, $\calG_X$ is the inverse image of the torsion subsheaf of
$\mathscr B_X/\calI_X$ under the quotient map from $\mathscr B_X$.
Consequently, $\mathscr B_X/\calG_X$ is torsion-free and hence
locally free on the smooth curve $X_1$. Thus $\calG_X$ is the
smallest subbundle of $\mathscr B_X$ containing the evaluation
image. Its generic fiber is the span of the global sections over
$k(X_1)$. Saturation does not change this span, but $\calG_X$ need
not be globally generated.

Raynaud's pairing \cite[\S4.1]{Raynaud1982} is the perfect 
alternating morphism
\[
b_X:\mathscr B_X\otimes_{\mathcal O_{X_1}}\mathscr B_X
       \longrightarrow \Omega^1_{X_1/k},
 \qquad b_X(du,dv)=c_X(u\,dv).
\]
It is well-defined: replacing
$u$ by $u+F_X^\#h$ changes the expression by $h\,c_X(dv)=0$.
The identity $2 b_X(du,du) = 2 c_X(u du) = c_X(d(u^2))=0$ shows that the pairing is alternating 
since $p \ne 2$.
The pairing is equivariant under automorphisms: if 
$\gamma_1$ is the Frobenius twist of an automorphism 
$\gamma$ of $X$, then $b_X(\gamma^*\alpha,\gamma^*\beta)
=\gamma_1^*b_X(\alpha,\beta)$.
This follows from the defining formula and the compatibility 
of Cartier with pullback.
Raynaud's pairing on the bundle $\mathscr B_X$ also leads to a pairing
of its global sections with values in $H^0(X_1,\Omega^1_{X_1/k})$;
see \cite[Proposition~1.2.1.1 and \S2.2]{Tong2010}.
At the generic point of $X_1$, choosing a nonzero rational
differential on $X_1$ identifies this pairing with a nondegenerate
alternating $k(X_1)$-bilinear form on the $(p-1)$-dimensional
generic fiber of $\mathscr B_X$.

\subsection{The hyperelliptic involution}

\begin{lemma}\label[lemma]{lem:hyperelliptic-isotropy}
If $X$ is hyperelliptic, then the pairing of any two global
sections of $\mathscr B_X$ is zero. The subbundle $\calG_X$ is
therefore isotropic for $b_X$.
\end{lemma}

\begin{proof}
Let $\iota, \iota_1$ be the hyperelliptic involutions on $X$ and $X_1$
respectively. Invariant regular differentials descend through
the tame quotient to $\bbP^1$, which has no nonzero regular
differentials. Therefore, the involutions act by $-1$ on 
$H^0(X,\Omega_X^1)$ and $H^0(X_1,\Omega^1_{X_1})$ respectively
(since $p\ne2$). 
For global sections $\alpha,\beta$ of $\mathscr B_X$, equivariance gives
\[
 \iota_1^*b_X(\alpha,\beta)
 =b_X(\iota^*\alpha,\iota^*\beta)
 =b_X(-\alpha,-\beta)=b_X(\alpha,\beta).
\]
The value is thus invariant, but as a regular differential on $X_1$
it is also anti-invariant. It follows that $b_X(\alpha,\beta)=0$.
By bilinearity the pairing vanishes on the generic span of the
global sections. Its restriction to $\calG_X\otimes\calG_X$
vanishes generically and hence vanishes as a morphism of vector
bundles. This proves isotropy.
\end{proof}

\begin{example}\label[example]{ex:quartic}
This is generally not true for non-hyperelliptic curves. For example, 
in characteristic $3$, on the smooth
plane quartic $X\colon x^4+y^4+1=0$, the regular exact differentials
$\alpha=dx/y^3=d(x/y^3)$ and
$\beta=x\,dx/y^3=d(x^2/(2y^3))$ satisfy
$b_X(\alpha,\beta)=c_X(x^2dx/y^6)=dx_1/y_1^2\ne0$.
\end{example}

\subsection{Eigenspaces for inverse characters}

Let $\pi:X\to\bbP^1$ be a connected tame cyclic cover of degree
$N$, and let $G=\Gal(X/\bbP^1)$. Fix a faithful character
$\chi:G\to k^\times$. For $0\le i<N$, define the eigenspace
\[
 D_i=\{\omega\in H^0(X,\Omega_X^1):
        \delta^*\omega=\chi(\delta)^i\omega
        \text{ for all }\delta\in G\},
 \qquad W_i=D_i\cap\ker\mathcal C_X.
\]
Choose a generator $\gamma$ of $G$ and set $\zeta_N=\chi(\gamma)$.
Then $\gamma^*$ acts on $D_i$ by $\zeta_N^i$. The invariant
eigenspace $D_0$ is zero. We regard $W_i$ as a subspace of
$H^0(X_1,\mathscr B_X)$, and let $\calG_i$ be the saturation of
the image of its evaluation map
$W_i\otimes_k\mathcal O_{X_1}\to\mathscr B_X$.
Write $\rho_i=\rank\calG_i$.

\begin{prop}\label[prop]{prop:inverse-ranks}
For $1\le i<N$, the subbundles $\calG_i$ and $\calG_{N-i}$ are
orthogonal for $b_X$, and
\begin{equation}\label{eq:inverse-ranks}
 \rho_i+\rho_{N-i}\le p-1.
\end{equation}
If $N$ is even, then $\calG_{N/2}$ is isotropic and
$\rho_{N/2}\le(p-1)/2$.
\end{prop}

\begin{proof}
For $\alpha\in W_i$ and $\beta\in W_{N-i}$, equivariance of the pairing gives
\[
 \gamma_1^*b_X(\alpha,\beta)
   =b_X(\zeta_N^i\alpha,\zeta_N^{-i}\beta)
   =b_X(\alpha,\beta),
\]
where $\gamma_1$ is the Frobenius twist of the automorphism $\gamma$ of $X$.
Thus $b_X(\alpha,\beta)$ is an invariant regular differential on $X_1$, and hence 
descends to a rational differential on $(\bbP^1)^{(p)}$. If we show that this 
differential is regular, then we will know that it vanishes, since $\bbP^1$ has 
no nonzero regular differentials.

At a point of $X_1$ with ramification index $e$, choose parameters
$u$ upstairs and $t$ downstairs in the completed local rings,
with $\pi_1^*t=u^e$. Tameness gives $p\nmid e$, and a generator
of inertia acts by $u\mapsto\zeta_eu$. If a differential
$\sum_{j\ge0}c_j u^jdu$ is invariant, then it forces $c_j=0$ unless 
$e\mid j+1$. Therefore the differential can be written as
\[
 u^{e-1}h(u^e)\,du
   =\pi_1^*\left(\frac{h(t)}e\,dt\right),
 \qquad h(t)\in k[[t]].
\]
showing that the differential downstairs is regular. Therefore, $b_X(\alpha,\beta)=0$.

The restriction of $b_X$ to $\calG_i\otimes\calG_{N-i}$ vanishes
on the generic fibers, and hence everywhere, just as in the proof
of \cref{lem:hyperelliptic-isotropy}. At the generic point,
$\calG_i\subseteq\calG_{N-i}^{\perp}$. Perfectness gives
$\rank\calG_{N-i}^{\perp}=p-1-\rho_{N-i}$, proving the
inequality in \eqref{eq:inverse-ranks}.
If $N$ is even, the character $\chi^{N/2}$ is its own inverse, so
$\calG_{N/2}$ is isotropic. Since $b_X$ is a perfect alternating
pairing on the rank-$(p-1)$ bundle $\mathscr B_X$, this isotropic
subbundle has rank at most $(p-1)/2$.
\end{proof}

\begin{remark}
For $N=2$, all regular differentials belong to $D_1$, so
$\calG_1=\calG_X$. Thus the even-degree case of
\cref{prop:inverse-ranks} gives the 
hyperelliptic isotropy of \cref{lem:hyperelliptic-isotropy}.
\end{remark}

\begin{remark}
The subbundles $\calG_i$ for distinct eigenspaces need not form a direct sum.
Although the eigenspaces $D_i$ are independent over $k$, their generic
spans are taken over $k(X_1)$, and rational-function coefficients
can change characters. Orthogonality gives a bound on the sum of
the ranks of the two subbundles $\calG_i$ and $\calG_{N-i}$; a half-rank bound 
requires the subbundle to be orthogonal to itself.
\end{remark}

\section{A description of the Cartier kernel in terms of polynomial syzygies}
\label{sec:polynomial}

We first give a polynomial counting argument for a finite separable
map $X\to\bbP^1$. It will be used here for hyperelliptic curves and,
in \cref{sec:cyclic}, separately for the eigenspaces of a tame cyclic
cover. After this common argument, we specialize to hyperelliptic curves 
and keep a fixed hyperelliptic model throughout the rest of the section.

Put $R=k[x]$, $K=k(x)$, and $K_0=k(x^p)$, the subfield of $p$-th
powers in $K$. The elements $1,x,\ldots,x^{p-1}$ form a basis of
$K$ over $K_0$. Since $k$ is perfect, every polynomial has a unique
expression
\begin{equation}\label{eq:Q}
 Q_{\vv}(x)=\sum_{j=0}^{p-1}v_j(x)^p x^j,
 \qquad \vv=(v_0,\ldots,v_{p-1})\in R^p.
\end{equation}
The same formula defines $Q_{\vv}\in K$ for $\vv\in K^p$.
For an integer $d\ge0$, write $R_{<d}=k[x]_{\le d-1}$.
The inequality $\deg Q_{\vv}<d$ is equivalent to
$\deg v_j\le\lfloor(d-1-j)/p\rfloor$ for every $j$.

\begin{prop}\label[prop]{prop:polynomial-count}
Let $x:X\to\bbP^1$ be a finite separable map, let $\eta$ be a
nonzero rational differential on $X$, and let $W$ be a $k$-subspace
of $R_{<d}\eta$. Define
\[
 \mathcal V_W=\{\vv\in R^p:Q_{\vv}\eta\in W\},
 \qquad s=\dim_{k(X)^p}\Span_{k(X)^p}W.
\]
Then $\dim_K\Span_K\mathcal V_W=s$ and
$\dim_k W\le\beta_p(d,s)$.
\end{prop}

\begin{proof}
The map $\vv\mapsto Q_{\vv}$ is additive and satisfies
$Q_{a\vv}=a^pQ_{\vv}$ for $a\in K$. Since Frobenius is an
isomorphism $K\to K_0$, it carries $K$-linear independence of
coefficient vectors to $K_0$-linear independence of their images.
The extension $k(X)^p/K_0$ is separable, whereas $K/K_0$ is purely
inseparable. They are linearly disjoint, so elements of $K$ that
are independent over $K_0$ remain independent over $k(X)^p$.
Multiplication by the nonzero differential $\eta$ preserves this
independence. This proves the rank identity. The correspondence
$\mathcal V_W\to W$ is Frobenius-semilinear and bijective, so
$\dim_k\mathcal V_W=\dim_kW$ as well.

Put $U=\Span_K\mathcal V_W$. If $s=0$, then $W=0$ and the
inequality is immediate. Otherwise, place a $K$-basis of $U$ in
the rows of an $s\times p$ matrix. A nonzero $s\times s$ minor
determines a set $I$ of $s$ column indices for which the coordinate
projection
\[
 \pi_I:U\longrightarrow K^s,\qquad
 (v_0,\ldots,v_{p-1})\longmapsto(v_j)_{j\in I}
\]
is injective. Its restriction to $\mathcal V_W$ is therefore an
injective $k$-linear map.

Write $d=pq+r$ with $0\le r<p$. The allowable polynomial space
for $v_j$ has dimension $q+1$ when $j<r$ and $q$ when $j\ge r$.
Consequently,
\[
 \dim_kW=\dim_k\mathcal V_W
 \le\sum_{j\in I}\bigl(q+\mathbf1_{j<r}\bigr)
 \le sq+\min\{s,r\}=\beta_p(d,s).
\]
\end{proof}

\subsection{Polynomial syzygies and the regularity bound}

For the remainder of this section, let $X$ be a hyperelliptic curve
with equation $y^2=f(x)$, where $f\in k[x]$ is squarefree of degree
$D\in\{2g+1,2g+2\}$. Put $n=(p-1)/2$ and write
\begin{equation}\label{eq:components}
 H:=f^n=\sum_{j=0}^{p-1}h_j(x)^p x^j,
 \qquad \mathbf h=(h_0,\ldots,h_{p-1}).
\end{equation}
Since $y^p=yH$, coefficient selection for Cartier gives (see \cite[\S3.1]{AchterHowe2019}):
\begin{equation}\label{eq:cartier-formula}
 \mathcal C_X\left(Q_{\vv}\frac{dx}{y}\right)
 = \frac{1}{y} \mathcal C_X\left(Q_{\vv}H dx\right)
 =\left(\sum_{i=0}^{p-1}v_i h_{p-1-i}\right)\frac{dx}{y},
\end{equation}
where the last equality comes from identifying the monomials in 
the product $Q_{\vv} H$ whose exponent is congruent to $p-1$ modulo $p$.

\begin{definition}\label[definition]{def:syzygy-lattice}
The \emph{rational syzygy space} of the Cartier components is the
$K$-vector space
\[
 V_f=\left\{\vv\in K^p:
             \sum_{i=0}^{p-1}v_i h_{p-1-i}=0\right\}.
\]
Its \emph{polynomial syzygy lattice} is the $R$-module
$L_f=V_f\cap R^p$.
\end{definition}

The row defining $V_f$ is nonzero, so the rational syzygy space in
\cref{def:syzygy-lattice} has dimension $p-1$.
The regular differentials on $X$ are $Q\,dx/y$ with $Q\in R_{<g}$.
Since the terms in \eqref{eq:Q} have disjoint exponent classes
modulo $p$, the condition $\deg Q_{\vv}\le g-1$ is equivalent to
\begin{equation}\label{eq:regularity-bound}
 \deg v_i\le\left\lfloor\frac{g-1-i}{p}\right\rfloor
 \qquad(0\le i<p).
\end{equation}
We call the inequalities in \eqref{eq:regularity-bound} the
\emph{regularity bound}.

\begin{definition}\label[definition]{def:regular-kernel}
Let $\calK_f$ be the $k$-vector space of polynomial syzygies
satisfying the regularity bound:
\[
 \calK_f=\left\{\vv\in L_f:
 \deg v_i\le\left\lfloor\frac{g-1-i}{p}\right\rfloor
                         \text{ for }0\le i<p\right\}.
\]
We denote its $K$-linear span in $V_f$ by $U_f=\Span_K\calK_f$.
\end{definition}

The map $\vv\mapsto Q_{\vv}\,dx/y$ identifies the space $\calK_f$
in \cref{def:regular-kernel}, Frobenius-semilinearly, with the
kernel of Cartier on regular differentials. Thus
\begin{equation}\label{eq:intrinsic-rank}
 a(X)=\dim_k\calK_f,\qquad \rho(X)=\dim_KU_f.
\end{equation}
For the second identity, \cref{prop:polynomial-count} applies
with $\eta=dx/y$ and $d=g$. More explicitly, its rank calculation is
\begin{align*}
\begin{split}
 \dim_KU_f
 &=\dim_{K_0}\Span_{K_0}\{Q_{\vv}:\vv\in\calK_f\}
  =\dim_{k(X)^p}\Span_{k(X)^p}\{Q_{\vv}:\vv\in\calK_f\}\\
 &=\dim_{k(X)^p}\Span_{k(X)^p}
                  \{Q_{\vv}\,dx/y:\vv\in\calK_f\}=\rho(X).
\end{split}
\end{align*}
The first equality uses the Frobenius-semilinear coefficient map;
the second uses the linear disjointness of $K$ and $k(X)^p$ over
$K_0$. The last dimension is the generic rank of the evaluation
image, because relative Frobenius identifies $k(X_1)$ with
$k(X)^p$. We are changing the scalar field on the polynomials
$Q_{\vv}$, not on the coefficient vectors $\vv$ themselves.
Polynomial multiples of a single vector can contribute several
dimensions to $a(X)$ while spanning only one direction in $U_f$.

\subsection{Centered heights and polynomial lattices}
\label{subsec:heights}

We want to count the polynomial syzygies satisfying the regularity
bound. A shift of the usual degree of each coordinate expresses
that bound as a single cutoff. The polynomial-lattice facts below
provide a basis adapted to this cutoff; they apply to any subspace
of $K^p$.

For $a\in K^\times$, let $\deg a$ be the degree of its numerator
minus that of its denominator. Define the centered shifts and the
\emph{centered Cartier length} by
\begin{equation}\label{eq:centered-height}
 t_i=\frac{i-n}{p},\qquad
 \htC(\vv)=\max_{0\le i<p}(\deg v_i+t_i),\qquad
 \htC(0)=-\infty.
\end{equation}
This length can be negative and satisfies
$\htC(a\vv)=\deg a+\htC(\vv)$. It is a function of the actual
vector: we do not divide a polynomial vector by its common
factors when imposing a length bound.

The height of a subspace uses a different normalization. Let
$V\subseteq K^p$ have dimension $d$, and let $M$ be a matrix whose
rows are a $K$-basis of $V$. For a subset
$I\subseteq\{0,\ldots,p-1\}$ of cardinality $d$, let $M_I$ be the
square submatrix in the columns indexed by $I$. The determinants
$\det M_I$ are its Pl\"ucker coordinates. A change of basis
multiplies all of them by the same element of $K^\times$.
After clearing denominators and dividing by the common polynomial
factor, we obtain a tuple $(P_I)_{|I|=d}$ of polynomials with
greatest common divisor $1$. It is unique up to multiplication by
an element of $k^\times$. These are the \emph{primitive polynomial
Pl\"ucker coordinates} of $V$.

Define the centered Cartier height of the subspace by
\begin{equation}\label{eq:subspace-height}
 \hC(V)=\max_{|I|=d}\left(\deg P_I+\sum_{i\in I}t_i\right),
 \qquad \hC(0)=0.
\end{equation}
The remaining ambiguity does not change any degree, so the height
is well-defined. In particular, the height of a line is the
centered Cartier length of a primitive polynomial generator,
not the length of an arbitrary vector on that line.

Let $L=V\cap R^p$. This is a free $R$-module of rank $d$.
It is saturated in $R^p$: if $a\vv\in L$ for $0\ne a\in R$ and
$\vv\in R^p$, then $\vv\in V$ and hence $\vv\in L$.
The quotient $R^p/L$ is therefore torsion-free, and thus free over
the principal ideal domain $R$. A basis of $L$ extends to a basis
of $R^p$, so its maximal minors are primitive.
For a real number $T$, put $L_{\le T}=\{\vv\in L:\htC(\vv)\le T\}$.
The $j$-th successive minimum $\mu_j(L)$ is the least $T$ for which
$L_{\le T}$ contains $j$ vectors independent over $K$.
These are shifted versions of standard polynomial-lattice
constructions \cite{Forney1975,Thunder1995,Thunder2008}.

\begin{prop}\label[prop]{prop:reduced-basis}
Let $V\subseteq K^p$ have dimension $d$ and set $L=V\cap R^p$.
\begin{enumerate}[(i)]
\item There is an $R$-basis $\mathbf b_1,\ldots,\mathbf b_d$ of $L$
whose centered Cartier lengths are the successive minima
$\mu_1\le \mu_2 \le \cdots\le\mu_d$, and for which
\begin{equation}\label{eq:predictable-degree}
 \htC\left(\sum_j a_j\mathbf b_j\right)
       =\max_{a_j\ne0}(\deg a_j+\mu_j)
\end{equation}
for all $a_1,\ldots,a_d\in R$, not all zero.
\item The subspace height is
$\hC(V)=\sum_{j=1}^d\mu_j$.
\item For every real $T$,
$\dim_kL_{\le T}=\sum_{j=1}^d\max\{0,1+\lfloor T-\mu_j\rfloor\}$.
\end{enumerate}
\end{prop}

\begin{proof}
For a nonzero polynomial vector $\mathbf b$, its \emph{leading
position} $\ell(\mathbf b)$ is the index attaining the maximum in
\eqref{eq:centered-height}. It is unique because the shifts $t_i$
are pairwise distinct modulo $\Z$.

\emph{(i)} Suppose two vectors $\mathbf b,\mathbf c$ of an
$R$-basis have the same leading position $i$, with
$\deg b_i\le\deg c_i$. Replacing $\mathbf c$ by
$\mathbf c-\lambda x^{\deg c_i-\deg b_i}\mathbf b$ for a suitable
$\lambda\in k^\times$ cancels its leading entry. This is a
unimodular basis operation and strictly decreases its length.
Repeating the operation terminates: the sum of the basis lengths
strictly decreases in $p^{-1}\Z$ and is bounded below by $-dn/p$.
We obtain a basis with distinct leading positions, which we order
by lengths $\lambda_1<\cdots<\lambda_d$.

For nonzero $a_j\in R$, the lengths of the vectors
$a_j\mathbf b_j$ have distinct fractional parts. The summand of
greatest length therefore cannot cancel, giving
\[
 \htC\left(\sum_j a_j\mathbf b_j\right)
       =\max_{a_j\ne0}(\deg a_j+\lambda_j).
\]
If a vector has length less than $\lambda_j$, its expansion in
this basis has zero coefficients at all indices $i\ge j$.
Such vectors span a $K$-space of dimension at most $j-1$.
On the other hand, $\mathbf b_1,\ldots,\mathbf b_j$ are independent
and have lengths at most $\lambda_j$. Hence $\mu_j=\lambda_j$,
which proves both assertions in (i).

\emph{(ii)} Write $i_j=\ell(\mathbf b_j)$ and consider the minor
on the columns $I=\{i_1,\ldots,i_d\}$. The term
$\prod_j b_{j,i_j}$ has degree $\sum_j\mu_j-\sum_{i\in I}t_i$.
Every other term uses a nonleading entry in at least one row;
after adding the same column shifts, its degree is strictly less
than $\sum_j\mu_j$. Thus this term cannot cancel in the determinant.
No maximal minor has degree plus its column shifts greater than
$\sum_j\mu_j$. Since the basis minors are primitive, they compute
\eqref{eq:subspace-height}, proving (ii).

\emph{(iii)} By \eqref{eq:predictable-degree},
$\sum_j a_j\mathbf b_j$ belongs to $L_{\le T}$ precisely when
$\deg a_j\le\lfloor T-\mu_j\rfloor$ for every nonzero $a_j$.
Consequently,
$\{x^e\mathbf b_j:1\le j\le d,\ 0\le e\le\lfloor T-\mu_j\rfloor\}$
is a $k$-basis of $L_{\le T}$, omitting indices with $T<\mu_j$.
Counting these vectors proves (iii).
\end{proof}

The centering also gives a useful duality. Define the
coefficient-reversal pairing and the corresponding orthogonal
complement by
\[
 \langle\vv,\mathbf w\rangle_{\mathrm{Car}}
    =\sum_{i=0}^{p-1}v_iw_{p-1-i},\qquad
 V^{\perpC}=\{\mathbf w\in K^p:
       \langle\vv,\mathbf w\rangle_{\mathrm{Car}}=0
                         \text{ for every }\vv\in V\}.
\]
Let $J$ be the $p\times p$ anti-diagonal permutation matrix.
Then $V^{\perpC}=J(V^\perp)$, where $V^\perp$ is the ordinary
dot-product orthogonal complement. If $(P_I)$ are primitive
Pl\"ucker coordinates for $V$, those for $V^\perp$, indexed by
$I^c$, can be chosen as $\pm P_I$. Reversal sends $I^c$ to
$\{p-1-j:j\in I^c\}$. Since $t_{p-1-j}=-t_j$ and
$\sum_jt_j=0$, we get
\begin{align}\label{eq:centered-duality}
 \hC(V^{\perpC})
 =\max_{|I|=d}\left(\deg P_I+\sum_{j\in I^c}t_{p-1-j}\right)
 =\max_{|I|=d}\left(\deg P_I-\sum_{j\in I^c}t_j\right)
  =\hC(V).
\end{align}
The signs do not affect the degrees. This proves self-duality for
the coefficient-reversal pairing.

Replacing the centered shifts $(i-n)/p$ by the uncentered shifts
$i/p$ increases the length of every nonzero vector and each
successive minimum by $n/p$, and increases the height of a
$d$-dimensional subspace by $dn/p$. Increasing a length cutoff
by the same amount leaves the admissible polynomial vectors
unchanged. The ordinary polynomial length is related by
$\max_i\deg v_i=\lfloor\htC(\vv)+n/p\rfloor$.

\subsection{The \texorpdfstring{$a$}{a}-number in terms of successive minima}
\label{subsec:a-minima}

We apply the preceding results to the polynomial syzygy lattice
$L_f$ of \cref{def:syzygy-lattice}. By \eqref{eq:Q},
\begin{equation}\label{eq:degree-Q}
 \deg Q_{\vv}=p\htC(\vv)+n\qquad(\vv\ne0).
\end{equation}
Thus the space $\calK_f$ of \cref{def:regular-kernel} is given by
\begin{equation}\label{eq:height-cutoff}
 \calK_f=(L_f)_{\le T_g},\qquad T_g=\frac{g-1-n}{p}.
\end{equation}
In particular, the regularity cutoff is inclusive.

The vector $\mathbf h$ of Cartier components in
\eqref{eq:components} is primitive: a common irreducible factor
$u$ would imply $u^p\mid f^n$, contrary to squarefreeness and
$n<p$. Since $\deg f^n=nD$, \eqref{eq:degree-Q} gives
$\htC(\mathbf h)=n(D-1)/p$. Moreover,
$V_f=(K\mathbf h)^{\perpC}$, so \eqref{eq:centered-duality} and
\cref{prop:reduced-basis}(ii) give
$\sum_{j=1}^{p-1}\mu_j(L_f)=\hC(V_f)=n(D-1)/p$.
The individual minima determine the regular kernel as follows.

\begin{prop}\label[prop]{prop:exact-a}
Let $X\colon y^2=f(x)$ be a smooth hyperelliptic curve, and let
$\mu_1\le\mu_2\le\cdots\le\mu_{p-1}$ be the centered successive minima
of its polynomial syzygy lattice $L_f$. With $T_g=(g-1-n)/p$,
\begin{equation}\label{eq:a-minima}
 a(X)=\sum_{j=1}^{p-1}\max\{0,1+\lfloor T_g-\mu_j\rfloor\},
 \qquad \rho(X)=\#\{j:\mu_j\le T_g\}.
\end{equation}
Moreover, $\mu_{n+1}>T_g$, so only the first $n$ terms in the sum
can be nonzero.
\end{prop}

\begin{proof}
By \eqref{eq:height-cutoff}, $\calK_f=(L_f)_{\le T_g}$.
Applying \cref{prop:reduced-basis}(iii) and then
\eqref{eq:intrinsic-rank} gives
\[
 a(X)=\dim_k\calK_f
     =\sum_{j=1}^{p-1}\max\{0,1+\lfloor T_g-\mu_j\rfloor\}.
\]
For the reduced basis $\mathbf b_1,\ldots,\mathbf b_{p-1}$ in
\cref{prop:reduced-basis}(i), the length formula gives
$U_f=\Span_K\{\mathbf b_j:\mu_j\le T_g\}$. These vectors are independent over $K$, so the rank identity in
\eqref{eq:intrinsic-rank} proves the formula for $\rho(X)$.
Finally, \cref{prop:inverse-ranks} with $N=2$ gives $\rho(X)\le n$;
hence at most $n$ minima lie below or at $T_g$.
\end{proof}

\subsection{The alternating form in polynomial coordinates}
\label{subsec:coordinate-pairing}

We now describe the pairing using the same hyperelliptic equation.
Its elementary ingredient is coefficient extraction on the base
$\bbP^1$. Every $u\in K$ has a unique expression
$u=\sum_{i=0}^{p-1}a_i x^i$ with $a_i\in K_0$.
Define the $K_0$-linear map $\tau:K\to K_0$ by
$\tau(u)=a_{p-1}$. Equivalently, if
$u=\sum_i b_i^p x^i$ with $b_i\in K$, then
$\tau(u)=b_{p-1}^p$, not $b_{p-1}$.

Differentiation is $K_0$-linear, since each element of $K_0$ has
zero derivative. It sends $\sum_i a_i x^i$ to
$\sum_{i=1}^{p-1}i a_i x^{i-1}$. In particular,
\[
 E:=\ker\tau=\bigoplus_{i=0}^{p-2}K_0x^i
          =\operatorname{im}(d/dx).
\]
Indeed, a primitive of $\sum_{i=0}^{p-2}a_i x^i$ is
$\sum_{i=0}^{p-2}a_i x^{i+1}/(i+1)$; all these denominators
are invertible. The same calculation shows that
$\ker(d/dx)=K_0$. Thus every $u\in E$ has a rational primitive,
and any two primitives differ by an element of $K_0$.

For $u,v\in E$, choose $U\in K$ with $U'=u$ and define
$\Omega(u,v)=\tau(Uv)\in K_0$. Replacing $U$ by $U+c$, with $c\in K_0$, changes the value by
$c\tau(v)=0$. Thus $\Omega$ is well-defined and $K_0$-bilinear.
It is alternating, because
$\tau(UU')=\tfrac12\tau((U^2)')=0$.
If $u$ belongs to its radical, then for every $W\in K$ we have
$0=\Omega(u,W')=\tau(UW')=-\tau(uW)$.
For $u\ne0$, taking $W=x^{p-1}/u$ contradicts this equality.
Hence $\Omega$ is nondegenerate.

For the rational syzygy space $V_f$ of \cref{def:syzygy-lattice},
multiplication by $H=f^n$ gives the Frobenius-semilinear
isomorphism $\Phi_f:V_f\to E$, $\Phi_f(\vv)=Q_{\vv}H$, since
\[
 \tau(Q_{\vv}H)=\left(\sum_i v_i h_{p-1-i}\right)^p.
\]
Here $\Phi_f(a\vv)=a^p\Phi_f(\vv)$ for $a\in K$.
To obtain a $K$-bilinear form on $V_f$, we must also transport
the output of $\Omega$ from $K_0$ to $K$ by inverse Frobenius.
Every element of $K_0=K^p$ has a unique $p$-th root in $K$, so
\begin{equation}\label{eq:Omega-f}
 \Omega_f(\vv,\mathbf w)
  =\Omega(\Phi_f(\vv),\Phi_f(\mathbf w))^{1/p}
  =\tau(U_{\vv}Q_{\mathbf w}H)^{1/p},
 \qquad U_{\vv}'=Q_{\vv}H.
\end{equation}
The root converts the scalar factor $a^pb^p$ into $ab$; hence
$\Omega_f(a\vv,b\mathbf w)=ab\Omega_f(\vv,\mathbf w)$.
The form is nondegenerate and alternating. It is different from
the symmetric coefficient-reversal pairing used for heights.

Let $\pi:X\to\bbP^1$ be the hyperelliptic map and
$\pi_1:X_1\to(\bbP^1)^{(p)}$ its twist. On function fields,
relative Frobenius induces
$F_{\bbP^1}^{\#}:k(x_1)\hookrightarrow k(x)$ with
$F_{\bbP^1}^{\#}(x_1)=x^p$.
For $u\in k(x)$, let $u^{[p]}(x_1)$ be obtained by raising its
coefficients to their $p$-th powers and replacing $x$ by $x_1$.
Thus $F_{\bbP^1}^{\#}u^{[p]}=u^p$, and relative Cartier satisfies
\begin{equation}\label{eq:relative-coefficient-formula}
 c_{\bbP^1}\left(\sum_{i=0}^{p-1}u_i(x)^p x^i\,dx\right)
        =u_{p-1}^{[p]}(x_1)\,dx_1.
\end{equation}
The superscript $[p]$ denotes the Frobenius twist, not a $p$-th root.

Put $\eta_{\vv}=Q_{\vv}H\,dx$ and $\omega_{\vv}=Q_{\vv}\,dx/y$.
Since $y^p=yH$, we have
$\omega_{\vv}=y^{-p}\pi^*\eta_{\vv}=d(U_{\vv}/y^p)$.

Using the definition of Raynaud's pairing, the linearity of
relative Cartier, and its compatibility with pullback, we obtain
\begin{align}\label{eq:pairing-comparison}
\begin{split}
 b_X(\omega_{\vv},\omega_{\mathbf w})
 &=c_X\left(\frac{U_{\vv}}{y^p}\,\omega_{\mathbf w}\right)
 =c_X\left(\frac1{y^{2p}}
               \pi^*(U_{\vv}Q_{\mathbf w}H\,dx)\right)
 =\frac1{y_1^2}\pi_1^*c_{\bbP^1}(U_{\vv}Q_{\mathbf w}H\,dx)\\
 &=\frac1{y_1^2}\pi_1^*b_{\bbP^1}(\eta_{\vv},\eta_{\mathbf w})
 =\frac1{y_1^2}\pi_1^*
       \left(\Omega_f(\vv,\mathbf w)^{[p]}(x_1)\,dx_1\right)
 =\pi_1^*\left(
       \frac{\Omega_f(\vv,\mathbf w)^{[p]}(x_1)}{f^{[p]}(x_1)}
       \,dx_1\right).
\end{split}
\end{align}
The third equality uses $F_X^\#y_1=y^p$, and the fifth uses
\eqref{eq:relative-coefficient-formula} and
$\tau(U_{\vv}Q_{\mathbf w}H)=\Omega_f(\vv,\mathbf w)^p$
from \eqref{eq:Omega-f}.
The last uses $y_1^2=f^{[p]}(x_1)$.
Pullback by the separable map $\pi_1$ is injective on rational
differentials. Therefore \eqref{eq:pairing-comparison} vanishes
exactly when $\Omega_f(\vv,\mathbf w)=0$. By
\cref{lem:hyperelliptic-isotropy}, the regular kernel span $U_f$
of \cref{def:regular-kernel} is isotropic for this elementary form.

\section{Upper bounds for tame cyclic covers}\label{sec:cyclic}

\subsection{Regular differentials in the eigenspaces}

We follow Elkin's description of regular differentials and Cartier
in \cite[\S\S2--3]{Elkin2011}. We include the regularity calculation
to make explicit which polynomial degrees enter the bound.

Fix a connected tame cyclic cover $\pi:X\to\bbP^1$ of degree $N$,
and choose the coordinate $x$ so that infinity is unramified.
Let $G=\Gal(X/\bbP^1)$. A Kummer equation for the cover is
\[
y^N=f(x)=\prod_\alpha(x-\alpha)^{a_\alpha},
 \qquad 0<a_\alpha<N,
\]
where $\alpha$ runs over the finite branch points.
The integer $D=\deg f$ is divisible by $N$.
Choose a generator $\gamma$ of $G$ with
$\gamma^*y=\zeta_N^{-1}y$, where $\zeta_N$ is a primitive $N$-th
root of unity. Let $\chi:G\to k^\times$ send $\gamma$ to
$\zeta_N$. The rational differential $dx/y^i$ is then a
$\chi^i$-eigenvector. It need not itself be regular.

As in \cref{sec:orthogonality}, let $D_i$ be the $\chi^i$-eigenspace
of regular differentials, and let $d_i=\dim_kD_i$.
Every rational differential with this character is
$q(x)\,dx/y^i$ for some $q(x)\in k(x)$.
At a point above a branch point $\alpha$, put
$e_\alpha=N/\gcd(N,a_\alpha)$. The orders of $x-\alpha$, $y$, and
$dx$ are respectively $e_\alpha$, $e_\alpha a_\alpha/N$, and
$e_\alpha-1$. It follows that regularity there is equivalent to
\[
 e_\alpha\operatorname{ord}_\alpha(q)+e_\alpha-1
       -\frac{i e_\alpha a_\alpha}{N}\ge0,
 \quad\text{or equivalently}\quad
 \operatorname{ord}_\alpha(q)\ge\left\lfloor\frac{i a_\alpha}{N}\right\rfloor.
\]
At the other finite points, $q$ must have no pole. Thus $q$ is a
polynomial divisible by
\begin{equation}\label{eq:character-polynomials}
 P_i(x)=\prod_\alpha(x-\alpha)^{\lfloor i a_\alpha/N\rfloor}.
\end{equation}
Write $q=QP_i$ and put $\eta_i=P_i\,dx/y^i$.
At each point over infinity, $x$ has pole order $1$, $y$ has pole
order $D/N$, and $dx$ has pole order $2$. Hence $Q\eta_i$ is
regular there exactly when $\deg Q+\deg P_i\le iD/N-2$.
Consequently
\begin{equation}\label{eq:character-space}
 D_i=R_{<d_i}\eta_i,\qquad
 d_i=\frac{iD}{N}-\deg P_i-1
     =\sum_\alpha\left\{\frac{i a_\alpha}{N}\right\}-1.
\end{equation}
Here $\{u\}=u-\lfloor u\rfloor$ denotes fractional part.
A basis is $\eta_i,x\eta_i,\ldots,x^{d_i-1}\eta_i$, with no vectors
when $d_i=0$. This is the basis in \cite[\S2, (2.4)]{Elkin2011}.
Since $D_0=0$, one has $g=\sum_{i=1}^{N-1}d_i$. If another
coordinate has a branch point at infinity, its branch exponent is
included in the fractional-part sum. Thus the dimensions depend
only on the ramification type, not on the positions of the branch
points.

We next recall the Cartier action on this basis. Let $\sigma(i)$ be the
integer in $\{1,\ldots,N-1\}$ congruent to $p^{-1}i$ modulo $N$,
and set
\[
 \epsilon_i=\frac{p\sigma(i)-i}{N},\qquad
 H_i=\frac{P_i f^{\epsilon_i}}{P_{\sigma(i)}^p}.
\]
Then $0\le\epsilon_i<p$. We compute the orders of $H_i$ as a
rational function on $\bbP^1$ at a branch point $\alpha$:
\begin{align}
 \operatorname{ord}_\alpha(H_i)
 =\left\lfloor\frac{ia_\alpha}{N}\right\rfloor
    +\epsilon_i a_\alpha
    -p\left\lfloor\frac{\sigma(i)a_\alpha}{N}\right\rfloor
 =\left\lfloor\frac{p\sigma(i)a_\alpha}{N}\right\rfloor
    -p\left\lfloor\frac{\sigma(i)a_\alpha}{N}\right\rfloor
 =\left\lfloor p\left\{\frac{\sigma(i)a_\alpha}{N}\right\}\right\rfloor.
\end{align}
The second equality uses $i+N\epsilon_i=p\sigma(i)$.
This order lies between $0$ and $p-1$. In particular, $H_i$ is a
polynomial. The identity $y^{p\sigma(i)}=y^i f^{\epsilon_i}$ gives
\[
 Q\eta_i=
 \left(\frac{P_{\sigma(i)}}{y^{\sigma(i)}}\right)^p
       \pi^*(QH_i\,dx).
\]
Using the semilinearity of Cartier and its compatibility with
pullback, we obtain
\begin{equation}\label{eq:cyclic-cartier}
 \mathcal C_X(Q\eta_i)
   =\frac{P_{\sigma(i)}}{y^{\sigma(i)}}
       \pi^*\mathcal C_{\bbP^1}(QH_i\,dx).
\end{equation}
Thus Cartier maps $D_i$ to $D_{\sigma(i)}$. Decomposing $Q$ and
$H_i$ according to exponents modulo $p$ expresses its kernel by
polynomial syzygies, just as in the hyperelliptic formula
\eqref{eq:cartier-formula}.

\subsection{Combining the rank restrictions}

Write $W_i=D_i\cap\ker\mathcal C_X$ and $w_i=\dim_kW_i$.
There are two restrictions on $w_i$. Elkin's lower bound for the
rank of Cartier on each eigenspace is
\begin{equation}\label{eq:elkin-character-estimate}
 \rank(\mathcal C_X|_{D_i})
     \ge\min\{2\lfloor d_i/p\rfloor,d_{\sigma(i)}\}
\end{equation}
\cite[Corollary~4.2]{Elkin2011}. Equivalently, $w_i\le e_i$, with
$e_i$ defined in \eqref{eq:elkin-caps}. On the other hand,
\cref{prop:polynomial-count} bounds $w_i$ in terms of the generic
rank $\rho_i$ of its evaluation image. The latter ranks are
constrained in pairs by \cref{prop:inverse-ranks}.

\begin{proof}[Proof of \cref{thm:cyclic-intro}]
By \eqref{eq:character-space}, $W_i$ is a subspace of
$R_{<d_i}\eta_i$. Applying \cref{prop:polynomial-count} gives
$w_i\le\beta_p(d_i,\rho_i)$. The inequality in
\eqref{eq:elkin-character-estimate} gives $w_i\le e_i$ as well.
The rank vector $(\rho_i)_{i=1}^{N-1}$ belongs to
$\mathscr R_{p,N}$ by \cref{prop:inverse-ranks}. Since Cartier
maps the eigenspaces $D_i$ to $D_{\sigma(i)}$, and $\sigma$ is a permutation, 
its kernel is $\bigoplus_{i=1}^{N-1}W_i$. Therefore
\[
 a(X)=\sum_{i=1}^{N-1}w_i
 \le\sum_{i=1}^{N-1}\min\{\beta_p(d_i,\rho_i),e_i\}
 \le\max_{\mathbf s\in\mathscr R_{p,N}}
          \sum_{i=1}^{N-1}\min\{\beta_p(d_i,s_i),e_i\}
 =\mathsf B_p(\mathbf d).
\]
\end{proof}

For each admissible vector $\mathbf s$, the $i$-th summand is at
most $e_i$ and at most $\beta_p(d_i,s_i)$. Thus the combined bound
$\mathsf B_p$ is bounded above by both the pairing bound $\mathsf P_p$
and Elkin's bound $\mathsf E_p$:
$\mathsf B_p(\mathbf d)\le\min\{\mathsf P_p(\mathbf d),\mathsf E_p(\mathbf d)\}$.
This comparison does not require the two separate bounds to be
ordered. The combined bound uses whichever restriction is stronger
on each eigenspace, subject to the joint rank constraint.

The combined bound $\mathsf B_p$ can be calculated one pair of
inverse characters at a time. To make it explicit, define $b_i(s)=\min\{\beta_p(d_i,s),e_i\}$ for $0\le s\le p-1$.
These functions are nondecreasing. For $i\ne N-i$, the rank
parameters $s_i,s_{N-i}$ occur only in their common constraint
$s_i+s_{N-i}\le p-1$. Increasing either variable cannot decrease
the objective, so the maximum is attained when
$s_i+s_{N-i}=p-1$. Adding the results therefore gives
\begin{equation}\label{eq:combined-pair-formula}
 \mathsf B_p(\mathbf d)=
 \sum_{i=1}^{\lfloor(N-1)/2\rfloor}
   \max_{0\le s<p}\{b_i(s)+b_{N-i}(p-1-s)\}
 +\begin{cases}
    b_{N/2}(n),&N\text{ even},\\
    0,&N\text{ odd}.
  \end{cases}
\end{equation}
If $N$ is even, the character $\chi^{N/2}$ is its own inverse.
The corresponding constraint is $2s_{N/2}\le p-1$, or equivalently
$s_{N/2}\le n$. Since its contribution is nondecreasing in
$s_{N/2}$, the maximum is attained at $s_{N/2}=n$.

\subsection{Comparison with Elkin's upper bound for the \texorpdfstring{$a$}{a}-number}
\label{subsec:comparison}

We first evaluate the contribution of two eigenspaces for inverse
characters to the pairing bound $\mathsf P_p(\mathbf d)$.
For $d,e\ge0$, put
\[
 \Psi_p(d,e)=\max_{0\le s<p}
        \{\beta_p(d,s)+\beta_p(e,p-1-s)\}.
\]
Write $d=pq+r$ and $e=pu+t$, with $0\le r,t<p$. Then
\begin{equation}\label{eq:psi-closed}
 \Psi_p(d,e)=
 \begin{cases}
  d-q,&q>u,\\
  e-u,&u>q,\\
  (p-1)q+\min\{p-1,r+t\},&q=u.
 \end{cases}
\end{equation}
Indeed, increasing $s$ by one changes the objective by
$q-u+\mathbf1_{s<r}-\mathbf1_{p-1-s\le t}$.
If $q>u$, this is nonnegative, and $s=p-1$ gives $d-q$.
The case $u>q$ is symmetric. If $q=u$, the linear part is
constant and the largest remainder contribution is
$\min\{p-1,r+t\}$. Summing these contributions over pairs of
inverse characters, and adding $\beta_p(d_{N/2},n)$ when $N$ is
even, computes the pairing bound $\mathsf P_p(\mathbf d)$.

\begin{proof}[Proof of \cref{thm:comparison-intro}]
\emph{(i)} Suppose $2\le N<p$. We prove the stronger comparison
\begin{equation}\label{eq:degree-comparison}
 \mathsf P_p(\mathbf d)
 \le g-2\sum_{i=1}^{N-1}\lfloor d_i/p\rfloor
 \le\mathsf E_p(\mathbf d).
\end{equation}
The second inequality follows from
$\min\{2\lfloor d_i/p\rfloor,d_{\sigma(i)}\}
 \le2\lfloor d_i/p\rfloor$ in Elkin's bound $\mathsf E_p$.
For the first, let $N_i=N/\gcd(N,i)$ be the order of $\chi^i$.
By \eqref{eq:character-space},
\[
 d_i+1=\sum_\alpha\left\{\frac{i a_\alpha}{N}\right\},
 \qquad
 d_{N-i}+1=\sum_\alpha
                  \left\{\frac{(N-i)a_\alpha}{N}\right\}.
\]
Every nonzero fractional part
$\theta_\alpha=\{ia_\alpha/N\}$ contributing to the first sum
belongs to $\{1/N_i,\ldots,(N_i-1)/N_i\}$. The corresponding
summand in the second sum is $1-\theta_\alpha$. If
$\theta_\alpha=0$, both summands are zero. Thus
$\theta_\alpha\le(N_i-1)(1-\theta_\alpha)$ for each nonzero
contribution. Summing and using $N_i\le N<p$ gives
\begin{equation}\label{eq:inverse-dim-comparison}
 d_i+1\le(N_i-1)(d_{N-i}+1)
       \le(p-2)(d_{N-i}+1).
\end{equation}

For $i\ne N-i$, write $d=d_i=pq+r$ and $e=d_{N-i}=pu+t$.
If $q>u$, the inequality in \eqref{eq:inverse-dim-comparison}
implies
\[
 q\le(p-2)u+
       \left\lfloor\frac{(p-2)(t+1)-r-1}{p}\right\rfloor
 \le(p-2)u+t.
\]
Using \eqref{eq:psi-closed}, we obtain
\[
 d+e-2q-2u-\Psi_p(d,e)=(p-2)u+t-q\ge0.
\]
The case $u>q$ follows by interchanging $d$ and $e$.
If $q=u$, the difference is
$(p-3)q+\max\{0,r+t-(p-1)\}\ge0$.
For a self-inverse character, writing its dimension as $d=pq+r$
gives
\[
 d-2q-\beta_p(d,n)
  =\frac{p-3}{2}q+\max\{0,r-n\}\ge0.
\]
Summing over all pairs of inverse characters and, when $N$ is even,
adding the contribution of the self-inverse character $\chi^{N/2}$
proves the first inequality in \eqref{eq:degree-comparison}.

\emph{(ii)} We use the restrictions imposed by ramification on the dimensions
of eigenspaces for inverse characters. For a vector
$\mathbf u=(u_1,\ldots,u_{N-1})$, put
$M(\mathbf u)=\frac12\sum_{i=1}^{N-1}\max(u_i,u_{N-i})$.
Thus $M(\mathbf d)$ counts the larger dimension in each pair once;
a self-inverse character contributes half its dimension.
The bounded remainder terms in \eqref{eq:psi-closed} give
\begin{equation}\label{eq:asymptotic-pair}
 \mathsf P_p(\mathbf d)
       =\frac{p-1}{p}M(\mathbf d)+O_{p,N}(1).
\end{equation}

We compare the contribution of each branch point to $M$ with its
total contribution to the eigenspace dimensions. The constant
\[
 c_N=
 \begin{cases}
  \displaystyle\frac34-\frac1{4N},&N\text{ odd},\\[3mm]
  \displaystyle\frac34-\frac1{4(N-1)},&N\text{ even}
 \end{cases}
\]
bounds this ratio at every branch point. In particular, $c_N<3/4$.
We will prove
\begin{equation}\label{eq:signature-imbalance}
 M(\mathbf d)\le c_Ng+(c_N-\tfrac12)(N-1).
\end{equation}
Set $z_i=d_i+1$. By \eqref{eq:character-space},
\[
 \mathbf z=\sum_\alpha\mathbf v_\alpha,
 \qquad
 \mathbf v_\alpha=
     \left(\left\{\frac{ia_\alpha}{N}\right\}\right)_{i=1}^{N-1}.
\]
If the inertia order at $\alpha$ is $e$, write
$a_\alpha/N=b/e$ with $\gcd(b,e)=1$. Among $i=0,\ldots,N-1$,
each residue modulo $e$ occurs $N/e$ times, and multiplication by
$b$ permutes these residues. Removing $i=0$ removes one zero.
Hence the nonzero entries of $\mathbf v_\alpha$ are
$1/e,\ldots,(e-1)/e$, each repeated $N/e$ times, and
$\sum_i v_{\alpha,i}=(N/e)(e-1)/2$.
If the entry indexed by $i$ is nonzero, the entry indexed by $N-i$
is $1$ minus that entry. If the entry indexed by $i$ is zero,
both entries are zero. Therefore
\begin{align}
      \frac{M(\mathbf v_\alpha)}{\sum_i v_{\alpha,i}}
      =
      \frac{\frac{N}{2e}\displaystyle\sum_{j=1}^{e-1}
            \max\left\{\frac je,\frac{e-j}{e}\right\}}
      {\frac Ne\displaystyle\sum_{j=1}^{e-1}\frac je}
      =
      \frac{\displaystyle\sum_{j=1}^{e-1}\max\{j,e-j\}}{e(e-1)}
      =
      \begin{cases}
      \displaystyle\frac34-\frac1{4e}, & e\text{ odd},\\[3mm]
      \displaystyle\frac34-\frac1{4(e-1)}, & e\text{ even}.
      \end{cases}
\end{align}
Since $e$ divides $N$, the ratio
$M(\mathbf v_\alpha)/\sum_i v_{\alpha,i}$
is at most $c_N$ in either case.

The function $M$ is subadditive, since
$\max(a+b,c+d)\le\max(a,c)+\max(b,d)$. Since 
$\mathbf z=\sum_\alpha\mathbf v_\alpha$, 
subadditivity of $M$ gives
\[
 M(\mathbf z)\le\sum_\alpha M(\mathbf v_\alpha)
 \le c_N\sum_\alpha\sum_i v_{\alpha,i}
 =c_N\sum_i z_i=c_N(g+N-1).
\]
Writing $\mathbf1=(1,\ldots,1)$, we obtain
\begin{align*}
 M(\mathbf d)=M(\mathbf z-\mathbf1)
             =M(\mathbf z)-\frac{N-1}{2}
 \le c_N(g+N-1)-\frac{N-1}{2}
    =c_Ng+\left(c_N-\frac12\right)(N-1),
\end{align*}
which proves the inequality in \eqref{eq:signature-imbalance}.

Elkin's bound satisfies $\mathsf E_p(\mathbf d)\ge(p-2)g/p$.
Combining this with \eqref{eq:asymptotic-pair} and the inequality in
\eqref{eq:signature-imbalance} gives
\begin{align}
 \mathsf E_p(\mathbf d)-\mathsf P_p(\mathbf d)
 \ge\left(\frac{p-2}{p}-\frac{p-1}{p}c_N\right)g-O_{p,N}(1)
 =\left[\frac{p-5}{4p}
       +\frac{p-1}{p}\left(\frac34-c_N\right)\right]g-O_{p,N}(1).
\end{align}
For $p\ge5$, the coefficient of $g$ is positive. The bounded term
depends only on $p,N$, so the pairing bound $\mathsf P_p$ is
strictly smaller than Elkin's bound $\mathsf E_p$ once $g$ is
sufficiently large, uniformly over the ramification types.
\end{proof}

In characteristic $3$ there is no such general ordering of the
two separate bounds. As we will see below, the family in
\cref{subsec:cyclic-sharpness} has Elkin's bound $\mathsf E_3$
strictly smaller than the pairing bound $\mathsf P_3$ in
arbitrarily large genus. Its $a$-number nevertheless attains the
combined bound $\mathsf B_3$, which is strictly smaller than both.

\section{Families with large \texorpdfstring{$a$}{a}-numbers}
\label{sec:families}

We first determine the regular Cartier kernel of the hyperelliptic
binomial family and prove the sharpness assertions of
\cref{thm:main}. We then give two infinite families attaining the
combined bound $\mathsf B_p$ for tame cyclic covers. The examples
illustrate the improvement from combining the pairing bound
$\mathsf P_p$ with Elkin's bound $\mathsf E_p$ in both possible
orderings of the separate estimates.

\subsection{The binomial family}\label{subsec:binomial}

Recall that $n=(p-1)/2$. Fix $1\le s<p$ and nonzero
$A,B\in k[T]$, and suppose that
\begin{equation}\label{eq:binomial}
 X\colon y^2=f(x),\qquad f(x)=A(x^p)+x^sB(x^p)
\end{equation}
has squarefree right-hand side of degree
$D\in\{2g+1,2g+2\}$. The binomial expansion is
\[
 f(x)^n=\sum_{j=0}^n\binom nj
            A(x^p)^{n-j}B(x^p)^j x^{sj}.
\]
The residues $sj\pmod p$, for $0\le j\le n$, are distinct, and
each summand is nonzero because $A,B\ne0$ and $n<p$.
Consequently, the set of residues that occur among the exponents
of $f^n$ is exactly
\[
 \Lambda_s=\{r\in\{0,\ldots,p-1\}:
                r\equiv sj\pmod p\text{ for some }0\le j\le n\}.
\]
Thus $|\Lambda_s|=n+1$, and precisely the $n$ Cartier components
whose indices are outside $\Lambda_s$ vanish.

Define the corresponding coordinate indices by
$I_s=\{p-1-r:0\le r<p,\ r\notin\Lambda_s\}$.
For every $i\in I_s$, the coordinate vector $\mathbf e_i$ belongs
to the polynomial syzygy lattice $L_f$ of
\cref{def:syzygy-lattice}. Its polynomial multiples satisfy the
regularity bound precisely up to degree
$\lfloor(g-1-i)/p\rfloor$. Hence, in every genus,
\[
a(X)\ge a_{p,g,s}:=
 \sum_{\substack{0\le t<p\\t\notin\Lambda_s}}
             \left\lfloor\frac{g+t}{p}\right\rfloor.
\]
The summand indexed by $t$ is the dimension of the allowed
polynomial space in coordinate $i=p-1-t$.

\begin{prop}\label[prop]{prop:binomial-exact}
Let $X$ be a smooth genus-$g$ member of \eqref{eq:binomial} with
$g\ge p$. Then
\begin{equation}\label{eq:binomial-exact}
 U_f=\bigoplus_{i\in I_s}K\mathbf e_i,
 \qquad
 \calK_f=\bigoplus_{i\in I_s}
       R_{\le\lfloor(g-1-i)/p\rfloor}\mathbf e_i.
\end{equation}
In particular, $\rho(X)=n$ and $a(X)=a_{p,g,s}$.
\end{prop}

\begin{proof}
When $g\ge p$, all $n$ constant vectors $\mathbf e_i$, $i\in I_s$,
satisfy the regularity bound. Their $K$-span is contained in $U_f$.
By \eqref{eq:intrinsic-rank} and \cref{prop:inverse-ranks} with
$N=2$, the dimension of $U_f$ is at most $n$. The first equality
in \eqref{eq:binomial-exact} follows. Every regular kernel vector
is therefore supported on $I_s$. Conversely, a polynomial vector
supported on $I_s$ automatically satisfies the syzygy equation, so
only its coordinate degree bounds remain. This proves the second
equality and the formulas for $\rho(X)$ and $a(X)$.
\end{proof}

For fixed $p$, \cref{prop:binomial-exact} determines the regular kernel in
all sufficiently large genera. It also describes the first $n$
centered successive minima. If $i_1<\cdots<i_n$ are the elements
of $I_s$, then $\mu_j(L_f)=(i_j-n)/p$ for $1\le j\le n$.
Indeed, all constant syzygies satisfy the regularity bound when
$g\ge p$, so the constant syzygy space is exactly the $k$-span of
these coordinate vectors. Every nonconstant polynomial vector has centered length
at least $(n+1)/p$, whereas every constant vector has length at
most $n/p$. The constant vectors therefore supply the first $n$
minima. For $s=1$, one has $I_1=\{0,\ldots,n-1\}$, giving the
smallest possible first $n$ minima of a sublattice of $R^p$.
In this case, if $g=pm+r$ with $0\le r<p$, the exact formula is
\begin{equation}\label{eq:numerical-sharpness}
 a(X)=nm+\min\{n,r\}.
\end{equation}

\subsection{Sharpness and codimension}\label{subsec:codimension}

For every $g\ge p$, we now choose smooth even-degree models in the
binomial family. Write $D=2g+2=pq+t$, with $0\le t<p$.
If $t>0$, take $s=t$ and allow $\deg A,\deg B\le q$.
If $t=0$, take any $1\le s<p$ and allow
$\deg A\le q$, $\deg B\le q-1$.
The resulting projective coefficient space has dimension
\begin{equation}\label{eq:parameter-dimension}
 \delta=\begin{cases}
          2q+1,&t>0,\\
          2q,&t=0.
         \end{cases}
\end{equation}
Its squarefree locus of exact degree $D$ is nonempty. An example is
\begin{equation}\label{eq:binomial-smooth-example}
 f(x)=1+x^s+x^D.
\end{equation}
If $t=0$, then $f'=s x^{s-1}$ and $f(0)=1$, so $f$ is squarefree.
If $t=s>0$, then $f'=s x^{s-1}(1+x^{D-s})$; at any nonzero root
of $f'$, the value of $f$ is $1$. This proves squarefreeness in
both cases.

\begin{prop}\label[prop]{prop:codimension}
Let $g\ge p$, write $2g+2=pq+t$ with $0\le t<p$, and choose
$s=t$ if $t>0$, or any $1\le s<p$ if $t=0$. Let $\delta$ be the
dimension in \eqref{eq:parameter-dimension}. We have the following
upper bound on the codimension in $\calH_g$ of the stratum of
hyperelliptic curves $X$ satisfying $a(X)\ge a_{p,g,s}$:
\begin{equation}\label{eq:codimension}
 \codim_{\calH_g}\calH_g^{\ge a_{p,g,s}}
       \le 2g+2-\delta.
\end{equation}
\end{prop}

\begin{proof}
The smooth equation family has dimension $\delta$. Its image in
$\calH_g$ lies in $\calH_g^{\ge a_{p,g,s}}$ by
\cref{prop:binomial-exact}. The moduli map has fibers of dimension
at most $3$, so its image has codimension at most $2g+2-\delta$
in $\calH_g$.
\end{proof}

\begin{proof}[Proof of \cref{thm:main}]
For a hyperelliptic curve, \cref{prop:inverse-ranks} with $N=2$
gives $\rho(X)\le n$. Apply \cref{prop:polynomial-count} to the
regular Cartier kernel, with $d=g$ and $\eta=dx/y$.
Using \eqref{eq:intrinsic-rank} and the fact that $\beta_p(g,s)$
is nondecreasing in $s$, we obtain
\[
 a(X)\le\beta_p(g,\rho(X))
      \le\beta_p(g,n)=nm+\min\{n,r\},\qquad g=pm+r.
\]
This proves both inequalities.

The equations in \eqref{eq:binomial-smooth-example} are defined over
$\F_p$ for every $g\ge p$, and \cref{prop:binomial-exact} gives
$\rho(X)=n$. If $g\equiv n\pmod p$, then $2g+2\equiv1\pmod p$
and we may take $s=1$. Equation \eqref{eq:numerical-sharpness}
gives equality in the bound for $a(X)$ as well. These genera
form an infinite progression, proving the sharpness assertions.
\end{proof}

\begin{remark}\label[remark]{rem:characteristic-two}
In characteristic $2$, the approach using locally exact differentials
still applies, but the argument needs Clifford's theorem in place
of isotropy with respect to Raynaud's pairing.
Indeed, $\mathscr B_X$ is a line bundle satisfying
$\mathscr B_X^{\otimes2}\simeq\Omega^1_{X_1}$
\cite[\S1.2.2]{Tong2010}, and $a(X)=h^0(X_1,\mathscr B_X)$.
Serre duality gives
\[
 h^1(X_1,\mathscr B_X)
 =h^0\!\left(X_1,\Omega^1_{X_1}\otimes\mathscr B_X^{-1}\right)
 =h^0(X_1,\mathscr B_X)=a(X).
\]
If $a(X)>0$, the line bundle $\mathscr B_X$ is effective and special.
Since its degree is $g-1$, Clifford's theorem
\cite[Chapter~IV, \S5 Theorem~5.4]{Hartshorne1977} gives
$a(X)\le\lfloor(g+1)/2\rfloor=\lceil g/2\rceil$. For hyperelliptic $X$, Elkin--Pries 
\cite[Theorem~1.3]{ElkinPries2013} in fact give an exact formula for $a(X)$ 
in terms of the ramification type. They also determine the isomorphism class 
of $\Jac(X)[2]$, and hence its Ekedahl--Oort type.
\end{remark}

\begin{proof}[Proof of \cref{thm:binomial-intro}]
The exact kernel formula follows from \cref{prop:binomial-exact}.
Since there are $n$ residues outside $\Lambda_s$, the defining sum
for $a_{p,g,s}$ is
$ng/p+O_p(1)$. Also
$\delta=2\lfloor(2g+2)/p\rfloor+O(1)=4g/p+O_p(1)$.
The inequality in \eqref{eq:codimension} therefore gives the stated
bound in $\calH_g$. Adding
$\codim_{\MM_g}\calH_g=g-2$ gives the bound in
\eqref{eq:intro-Mg} for $\MM_g$.
\end{proof}

\subsection{Infinite families attaining the bound for tame cyclic covers}
\label{subsec:cyclic-sharpness}

The following two families attain the combined bound
$\mathsf B_p$ in infinitely many genera. In characteristic $3$,
Elkin's bound $\mathsf E_3$ is smaller than the pairing bound
$\mathsf P_3$; in characteristic $5$, the ordering is reversed.
In both cases, the combined bound is strictly smaller than either
separate bound.

\begin{prop}\label[prop]{prop:sharp-families}
For every integer $t\ge1$, the smooth projective models in the
following table are nonhyperelliptic and have the indicated genus
and $a$-number, which equals the combined bound $\mathsf B_p(\mathbf d)$:
\begin{center}
\begin{tabular}{c c c c}
\toprule
$p$ & Equation for $X_t$ & $g(X_t)$ & $a(X_t)=\mathsf B_p(\mathbf d)$\\
\midrule
$3$ & $y^5=x^{15t+5}-1$ & $30t+6$ & $11t+3$\\
$5$ & $y^3=x^{15t+3}-1$ & $15t+1$ & $7t+1$\\
\bottomrule
\end{tabular}
\end{center}
The pairing bound $\mathsf P_p(\mathbf d)$ and Elkin's bound
$\mathsf E_p(\mathbf d)$ are respectively $14t+4$ and $13t+4$
for $p=3$, and $8t+1$ and $9t+1$ for $p=5$.
\end{prop}

\begin{proof}
Write each equation as $y^N=x^L-1$. In both cases $p\nmid NL$
and $N\mid L$. Thus the right-hand side is squarefree, the cover
is tame and connected, and infinity is unramified. In
\eqref{eq:character-polynomials} we have $P_i=1$, so
\[
 d_i=\frac{iL}{N}-1,\qquad
 D_i=\Span_k\left\{x^j\frac{dx}{y^i}:0\le j<d_i\right\}.
\]
By \eqref{eq:cyclic-cartier}, Cartier on a basis vector is obtained
by selecting terms of $x^j(x^L-1)^{\epsilon_i}$ whose exponents
are congruent to $p-1$ modulo $p$. Since $0\le\epsilon_i<p$,
all binomial coefficients in the expansion are nonzero. A term
contributes precisely when
\begin{equation}\label{eq:monomial-selection}
 j+\ell L\equiv p-1\pmod p,\qquad 0\le\ell\le\epsilon_i.
\end{equation}
For a fixed $j$ there is at most one such $\ell$, since $L$ is
prime to $p$. Distinct nonzero outputs are also independent:
a collision would give $j-j'=(\ell'-\ell)L$, impossible for
$j\ne j'$ because $0\le j,j'<d_i<L$. Hence the kernel has a
basis consisting of the individual basis differentials killed
by Cartier.

\Cref{tab:sharp-family-kernels} records the zero-column residues
from \eqref{eq:monomial-selection}. Counting their occurrences in
$0\le j<d_i$ gives $w_i=\dim_k(D_i\cap\ker\mathcal C_X)$.
\begin{table}[H]
\centering
\small
\begin{tabular}{c c c c c c c c}
\toprule
$p$&$N$&$i$&$d_i$&$\sigma(i)$&$\epsilon_i$
 &Kernel residues $j\bmod p$&$w_i$\\
\midrule
\multirow{4}{*}{$3$}&\multirow{4}{*}{$5$}
 &$1$&$3t$&$2$&$1$&$\{1\}$&$t$\\
 &&$2$&$6t+1$&$4$&$2$&$\varnothing$&$0$\\
 &&$3$&$9t+2$&$1$&$0$&$\{0,1\}$&$6t+2$\\
 &&$4$&$12t+3$&$3$&$1$&$\{1\}$&$4t+1$\\
\midrule
\multirow{2}{*}{$5$}&\multirow{2}{*}{$3$}
 &$1$&$5t$&$2$&$3$&$\{2\}$&$t$\\
 &&$2$&$10t+1$&$1$&$1$&$\{0,2,3\}$&$6t+1$\\
\bottomrule
\end{tabular}
\captionsetup{labelfont=up}
\caption{Cartier kernels in the $\chi^i$-eigenspaces for the families
in \cref{prop:sharp-families}. The kernel residues are the classes
of $j$ modulo $p$ for which $x^j\,dx/y^i$, with $0\le j<d_i$,
is killed by Cartier. The last column gives the kernel dimension.}
\label{tab:sharp-family-kernels}
\end{table}
Summing the $d_i$ and $w_i$ gives the genera and the values of
$a(X_t)$ in the statement.

We compute the combined bound $\mathsf B_p(\mathbf d)$ using
\eqref{eq:combined-pair-formula}. For $p=3$, the individual
kernel bounds are
$(e_1,e_2,e_3,e_4)=(t,2t+1,6t+2,4t+1)$.
For the inverse characters $\chi$ and $\chi^4$, the choices
$s=0,1,2$ give $4t+1,5t+1,t$; their maximum is $5t+1$.
For the inverse characters $\chi^2$ and $\chi^3$, the same choices
give $6t+2,5t+2,2t+1$, with maximum $6t+2$.
Thus the combined bound is $\mathsf B_3(\mathbf d)=11t+3$.

For $p=5$, the individual bounds are $(e_1,e_2)=(3t,6t+1)$.
For the inverse characters $\chi$ and $\chi^2$, the five choices
$s=0,\ldots,4$ give $6t+1$, $7t+1$, $6t+1$, $5t+1$, and $3t$,
respectively. Their maximum is $7t+1$, proving that the combined bound is
$\mathsf B_5(\mathbf d)=7t+1$.
In each case, the stated value of the pairing bound $\mathsf P_p$
is obtained from \eqref{eq:psi-closed}, and that of Elkin's bound
$\mathsf E_p$ is obtained by summing the $e_i$.

Finally, the differentials $dx/y^2$, $x\,dx/y^2$, and $dx/y$ are
regular in both families. Their ratios recover $x$ and $y$, so
the canonical map is birational. The curves are therefore
nonhyperelliptic.
\end{proof}

\printbibliography
\end{document}